\documentclass[12pt]{amsart}
\usepackage{amscd}
\usepackage{verbatim}
\usepackage{amssymb, amsmath, amsthm, amscd,ifthen}
\usepackage[cp866]{inputenc}

\usepackage{tikz}

\usepackage{amsmath,amssymb,amscd,}
\usepackage[all]{xy}

\usepackage{amssymb}

\newtheorem{thm}{Theorem}[section]

\newtheorem{prop}[thm]{Proposition}

\newtheorem{Ex}[thm]{Example}

\newtheorem{lemma}[thm]{Lemma}
\newtheorem{cor}[thm]{Corollary}

\theoremstyle{definition}

\newtheorem{dfn}[thm]{Definition}
\newtheorem{conj}{Conjecture}

\def\Z{\mathbb Z}

\def\comment#1{}
\begin{document}
\title{Swap action on   moduli spaces of  polygonal linkages }
\author{M. Khristoforov, G. Panina}
\begin{abstract} The basic object of the paper are the moduli spaces $M_{2}(L)$ and  $M_{3}(L)$ of a closed
polygonal linkage either in $\mathbb{R}^2$ or in $\mathbb{R}^3$. As
was originally suggested by G. Khimshiashvili, the space  $M_{2}(L)$
is equipped with the oriented area function $A$. In turn, we equip
the space $M_{3}(L)$ with the vector area function $S$.
 The latter are
generically  Morse functions, whose critical points  have a nice
description. In the paper, we define a \textit{swap action} (that
is, the action of some group generated by edge transpositions) on
the spaces $M_{2}(L)$ and  $M_{3}(L)$ which preserves the functions
$A$ and $S$ and the Morse points. We prove that the commutant of the
group acts trivially, present some computer experiments and
formulate a conjecture.
\end{abstract}
\maketitle
\section{Introduction}

We study the moduli spaces $M_{2}(L)$ and  $M_{3}(L)$ of a closed
polygonal linkage either in $\mathbb{R}^2$ or in $\mathbb{R}^3$.
These spaces attract special attention firstly because of practical
applications, and secondly because they can be equipped by
additional structures. In this respect we briefly mention the papers
by A. Klyachko \cite{klya}, and by M. Kapovich, J. Millson
\cite{KM}.

 In the paper, we consider the \textit{oriented
area function} $A$ defined on the space $M_{2}(L)$, and the\textit{
vector area function} $S$ defined on the space $M_{3}(L)$.
Generically,  $A$ is
 a Morse functions, and $S$ is a Morse-Bott function, whose critical points  have a nice
description. In the paper, we enrich this structure by defining a
\textit{swap action} (that is, the action of some group generated by
edge transpositions) on the  spaces $M_{2}(L)$ and  $M_{3}(L)$ which
preserves the functions $A$ and $S$ and the Morse points. We show
that this action factors through  a factor group of  \textit{the
group of pure balanced annular braids}. Besides, we prove that
commutant of the group acts trivially, present some computer
experiments and formulate a natural conjecture.


\textbf{Acknowledgments.} The first author was supported by the
Chebyshev Laboratory (Department of Mathematics and Mechanics,
St.-Petersburg State University) under RF government grant
11.G34.31.0026.

\section{Moduli space and oriented area}

A polygonal  \textit{ $n$-linkage} is a sequence of positive numbers
$l_1,\dots ,l_n$. It should be interpreted as a collection of rigid
bars of lengths $l_i$ joined consecutively by revolving joints in a
closed  chain. We study its flexes  with allowed self-intersections.
This is formalized in the following definition:
\begin{dfn}
For a  linkage $L$, \textit{a configuration} in the Euclidean space
$ \mathbb{R}^d$ is a sequence of points $R=(p_1,\dots,p_{n}), \ p_i
\in \mathbb{R}^d$ with $l_i=|p_i,p_{i+1}|,  \  l_n=|p_np_1|$.

 The  \textit{the moduli space of  }$L$ is the set $M_d (L)$ of all such
configurations modulo the action of orientation preserving
isometries.
\end{dfn}

In the paper we make use of the signed area function as the Morse
function on $M_{2}(L)$ and of the vector area function on
$M_{3}(L)$.

\bigskip

We start with the \textbf{2D}-case.

\begin{dfn} \label{Dfn_areaR2} The \textit{signed area} of a polygon $P \subset \mathbb{R}^2$ with the vertices \newline $p_i = (x_i,
y_i)$  is defined by
$$2A(P) = (x_1y_2 - x_2y_1) + \ldots + (x_ny_1 - x_1y_n).$$
\end{dfn}

\begin{dfn}
    A
polygon  $P$  is called \textit{cyclic} if all its vertices $p_i$
lie on a circle.
\end{dfn}

 Cyclic polygons   arise in the framework of
 the paper as critical points of the signed area:

\begin{thm}\label{Thm_crirical_are_cyclic}\cite{khipan}
  Generically, a polygon $P$ is a critical point of the
signed area function $A$  iff $P$ is a cyclic configuration.
         \qed
\end{thm}

\begin{figure}
\centering
\includegraphics[width=8 cm]{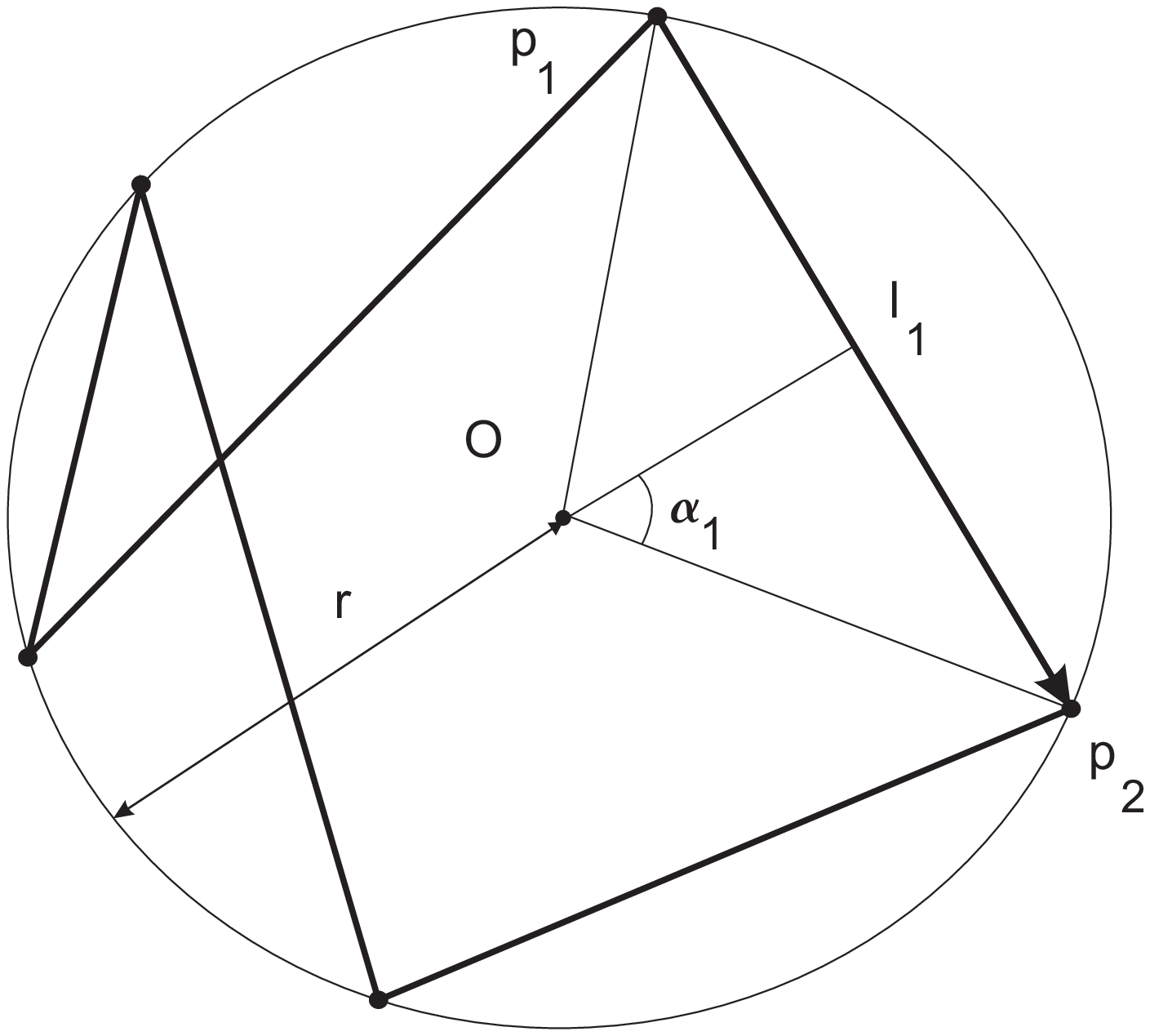}
\caption{Basic notation for a pentagonal cyclic configuration with
$E=(-1,-1,-1,1,-1)$}\label{Figure_notation}
\end{figure}

 The following notation (see Fig.\ref{Figure_notation}) is used throughout the paper for closed cyclic configurations:

$r=r(P)$ is the radius of the circumscribed circle.

A cyclic configuration  is called \textit{central} if one of its
edges contains $O$.

For a non-central configuration,  $\varepsilon_i$ is the
\textit{orientation} of the edge $p_ip_{i+1}$:

 $\varepsilon_i=\left\{
                       \begin{array}{ll}
                         1, & \hbox{if the center $O$ lies to the left of } p_ip_{i+1};\\
                         -1, & \hbox{if the center $O$ lies to the right of } p_ip_{i+1}.
                       \end{array}
                     \right.$

$E=E(P)=(\varepsilon_1,\dots,\varepsilon_n)$ is the string of
orientations of all the edges.

\newpage

Now we pass to the \textbf{3D}-case.

\begin{dfn} \label{Dfn_areaR3} The \textit{vector area} of a polygon $P \subset \mathbb{R}^3$ with the vertices \newline $p_i = (x_i,
y_i)$  is defined by
$$2\overrightarrow{S(P)}= p_1 \times p_2 + p_2 \times p_3+ \dots + p_n\times p_1,$$
$$2S(P)= |p_1 \times p_2 + p_2 \times p_3+ \dots + p_n\times p_1|.$$

\end{dfn}

\begin{figure}
\centering
\includegraphics[width=12 cm]{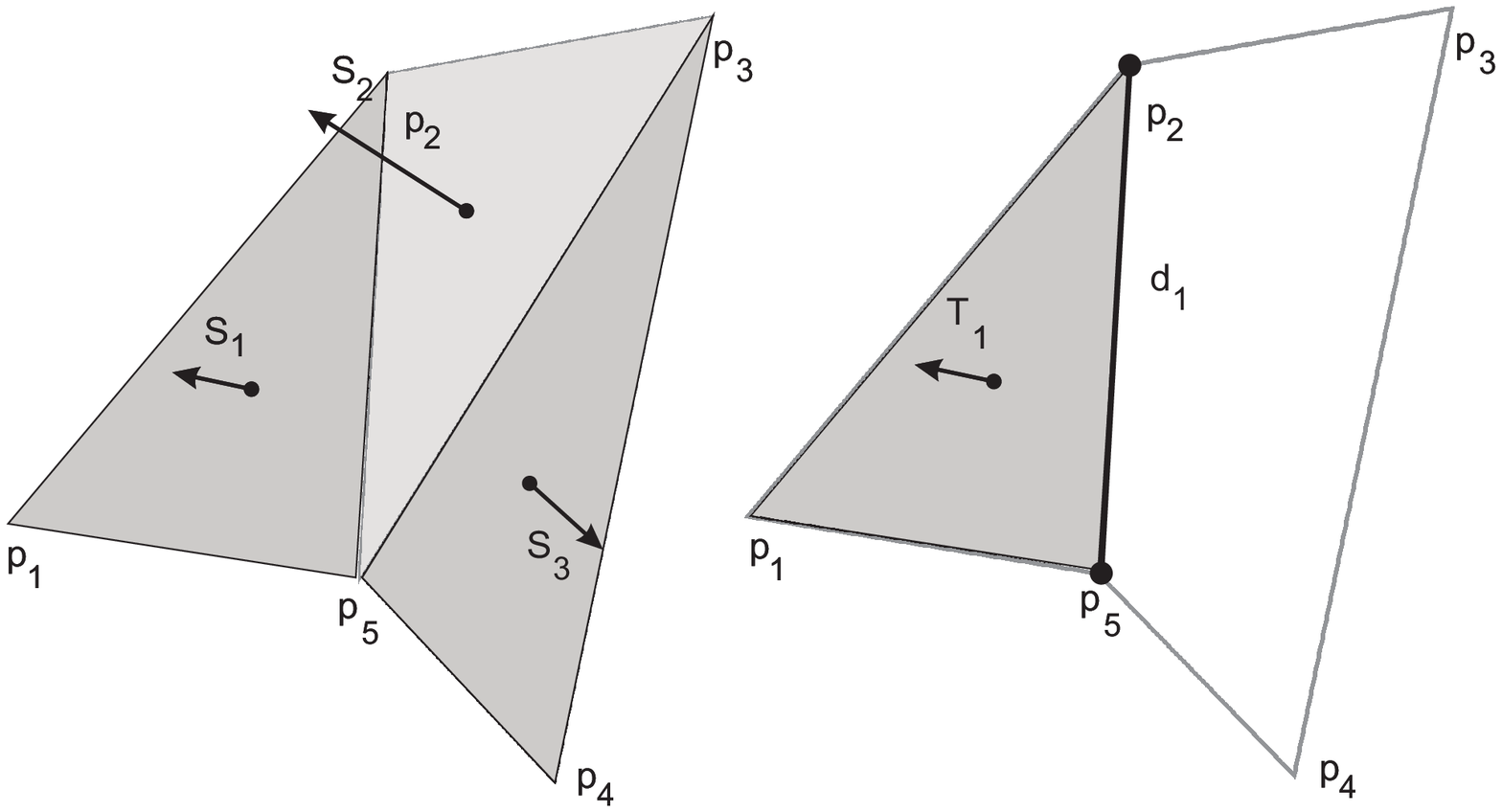}
\caption{}\label{Figure_3Darea}
\end{figure}

\begin{thm}\label{Thm_crirical_3D}
 Assume that $S(P)\neq 0$ for a configuration $P\in M_3(L)$. Generically,  $P$ is a critical point of the
vector area function $S$  if and only if  the two following
conditions hold:
\begin{enumerate}
    \item The orthogonal projection of $P$ onto the plane $\overrightarrow{S(P)}^\perp$ is a
    cyclic polygon.
    \item For every $i$, the vectors $\overrightarrow{T_i}$, $\overrightarrow{S}$, and $\overrightarrow{d_i}$ are
    coplanar.
\end{enumerate}
Here $\overrightarrow{d_i}$ is the $i$-th short diagonal,
$\overrightarrow{T_i}$ is the vector area of the triangle
$p_{i-1}p_ip_{i+1}$, see Fig. \ref{Figure_3Darea}, right.
\end{thm}
Proof. We list $(2n-6)$ flexes that generate  $(2n-6)$ elements of
the tangent space $T_P(M_3(L))$. Generically, these vectors are
linearly independent.  Therefore, the point $P\in M_3(L)$ is
critical if  and only if the function $S$ has a zero derivative in
all these directions.\begin{enumerate}
    \item Denote by $pr\ P$ the orthogonal projection of $P$ onto
    $\overrightarrow{S(P)}^\perp$. Each flex of $pr\ P$ in the plane $\overrightarrow{S(P)}^\perp$
    generates a flex of $P$ in the space $\mathbb{R}^3$. During the flex, we
    maintain the slopes of the edges with respect to the  plane
    $\overrightarrow{S(P)}^\perp$. Since $dim \ M_2(pr\ P)=n-3$, we can choose  $(n-3)$
    linearly independent tangent vectors of this type.
    \item Let us  bend   the triangle $T_i$ around  the diagonal
    $d_i$ keeping the rest of configuration $P$ frozen.  We  choose  $(n-3)$
    linearly independent tangent vectors of this type.
\end{enumerate}
The flexes of the first (respectively, second) type  provide the
statement 1 (respectively, statement 2) of the theorem. \qed

\section{Swap action}
\begin{figure}
\centering
\includegraphics[width=12 cm]{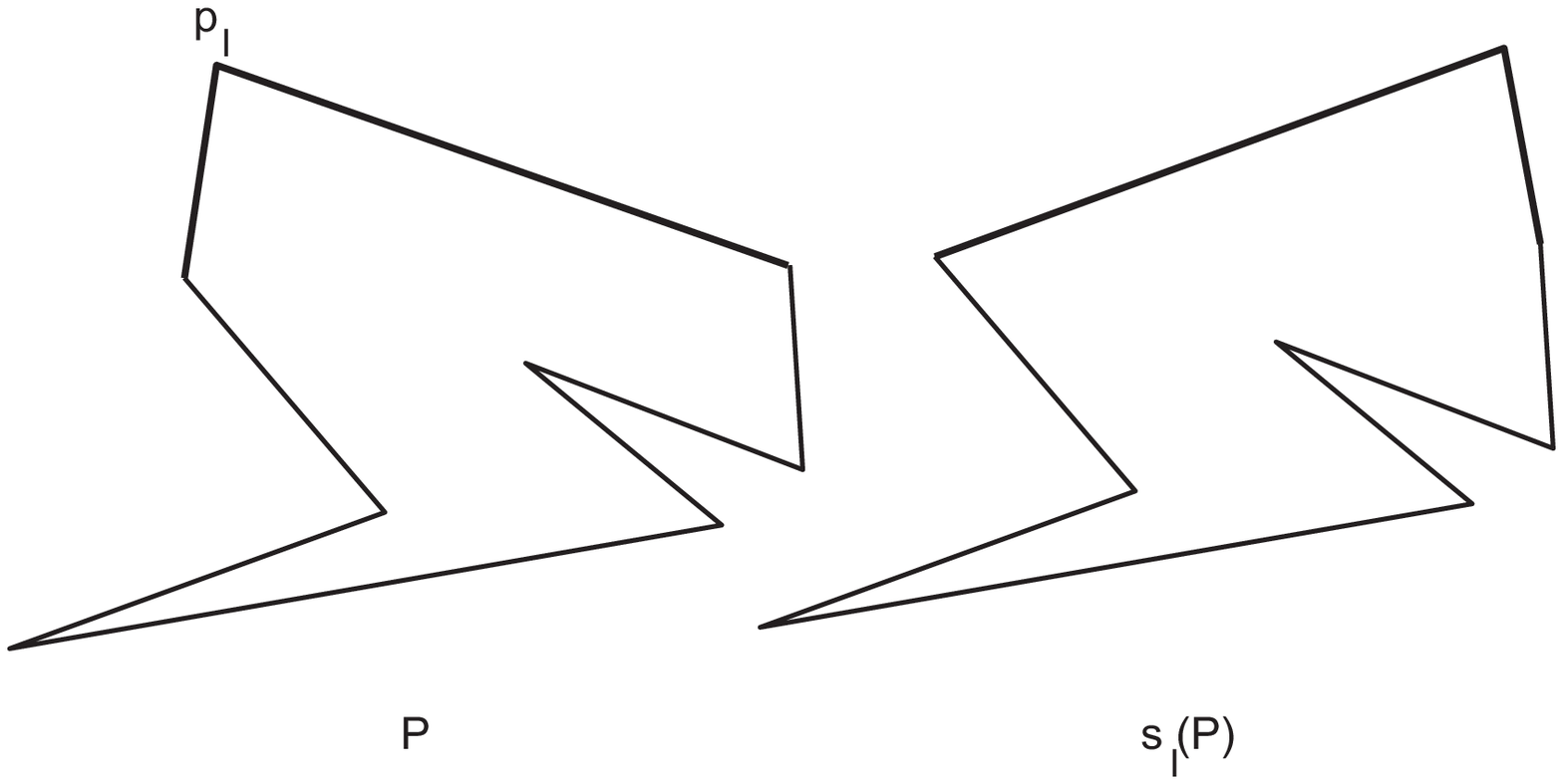}
\caption{}\label{Figure_swap}
\end{figure}

We assume that a polygonal linkage $L$ with $n$ edges and with all
$l_i$ different is fixed. We make a convention that the numbering is
modulo $n$, that is, for instance, $n+1=1$.
\begin{dfn}
Let $P\in M_{2,3}(L)$ be a polygon. For $i=1, \dots, n$, denote by
$s_i(P)$ the polygon obtained from $P$ by transposing of the two
edges adjacent to the vertex $p_i$ (see Fig. \ref{Figure_swap}).
For the dimension three, we assume that the new pair of edges lies
in the plane spanned by the old one.
\end{dfn}

We get a homeomorphism
$$s_i:M_{2,3}(L)\rightarrow M_{2,3}(\sigma_i L), $$
where the element of the symmetric group $\sigma_i \in S_n$ is a
transposition induced by $s_i$. Define by $F_n$ the free group whose
generators are the abstract symbols $s_i$.

\begin{lemma} \begin{enumerate}
    \item  The action of $F_n$ respects the  functions $A$ and $\overrightarrow{S}$.
    \item For $n=4$, the action of $F_4$ respects the volume of the
convex hull $V(Conv(P))$. \qed
\end{enumerate}
\end{lemma}

However, $F_n$ acts on the disjoint union of moduli spaces
$\bigsqcup M_{2,3}(\sigma_iL)$. We wish to restrict ourselves by
 just one moduli space. This means that we take
only those elements that take a configuration to the same moduli
space. We formalize this as follows: \comment{ We formalize this as
follows: Denote by $C_n\subset S_n$ the subgroup of the symmetric
group generated by the cycle $(2,3, \cdots,n,1)$.
There is a natural mapping to the factor group
$$\pi:SW_n\rightarrow S_n/C_n.$$
We are interested in the action of its kernel $SW_n^0$ on the moduli
space $ M_{2,3}( L).$ }
 There is a natural mapping to the
symmetrical group
$$\pi:F_n\rightarrow S_n,$$
which maps $s_i$ to $\sigma_i$. Clearly its kernel $F_n^0$ acts  on
the moduli space $ M_{2,3}( L).$

\begin{lemma} \label{Lemma4trivial}For a $4$-linkage $L$, the group $F_4^0$  acts trivially
on $M_{2,3}(L)$.
\end{lemma}
Proof. (2D). For a $4$-gon  $P = (p_1, p_2, p_3, p_4)$ denote by
 $O=O( P)$ the  intersection point of perpendicular
bisectors to the segments $p_1p_3$ and $p_2p_4$. Denote also
$$r_i( P) = |Op_i|, \quad \beta_i( P) = \angle p_iOp_{i+1}.$$
The lemma  follows from the three geometrical observations:
\begin{enumerate}
\item A $4$-gon is completely defined by
$$r(P) = ((r_1(P), r_2(P), r_3(P), r_4(P)),\hbox{ and }
\beta(P) = (\beta_1(P),\beta_2(P),\beta_3(P),\beta_4(P)). $$
\item The action of $ F_n$ preserves the point  $O(P)$ and the vector $r(P)$.
\item The group  $ F_n$  acts on $\beta(P)$ by permutations:
$ \beta(s(P)) = \pi(s) \beta(P).$
\end{enumerate}

\bigskip

(3D).   By  analyticity reasons it is enough to prove that $
s=(s_1s_2)^3$ acts trivially on some  open subset $U$ of the space
of all $4$-gons.

 Take an equilateral $4$-gon  $P_0$ (that is,  a rhombus but not a
 square). The swap  $s$ obviously takes $P_0$ to itself. Now, let $P$
 be a quadrilateral close to $P_0$. Its image $sP$ is  close to $P$ and has the
 same values of  $S(P)$ and $V(Conv(P))$. By  continuity reasons, $sP=P$.
In other words, $ s=(s_1s_2)^3$ acts trivially on a neighborhood of
$P$ which is an open set. \qed

\bigskip

\begin{dfn}Denote by $ Stab = Stab (M_{2,3}(L))\subset F^0_n$  the pointwise stabilizer of
the space $M_{2,3}(L)$, that is, the the group of all elements with
the trivial action. Denote also the factor $F^0_n/Stab$ by
$SW_n=SW_n(L)$.
\end{dfn}

\begin{prop} Generically, the group $Stab$ does not depend on
$L$.\qed
\end{prop}

\begin{dfn}
 Define $R \subset F^0_n$ as the subgroup generated by the elements of the following three types:
\begin{enumerate}
    \item $s_i^2$,
    \item $s_is_js_i^{-1}s_j^{-1}$, whenever $|i-j|>1$, and
    \item $s_is_{i+1}s_is_{i+1}^{-1}s_i^{-1}s_{i+1}^{-1}$.
\end{enumerate}
\end{dfn}

\begin{prop} The group $R$ is a subgroup of the stabilizer
 $ Stab $. 
\end{prop}

Proof. The first two items are obvious. The third one follows from
Lemma \ref{Lemma4trivial}. \qed

\begin{figure}
\centering
\includegraphics[width=14 cm]{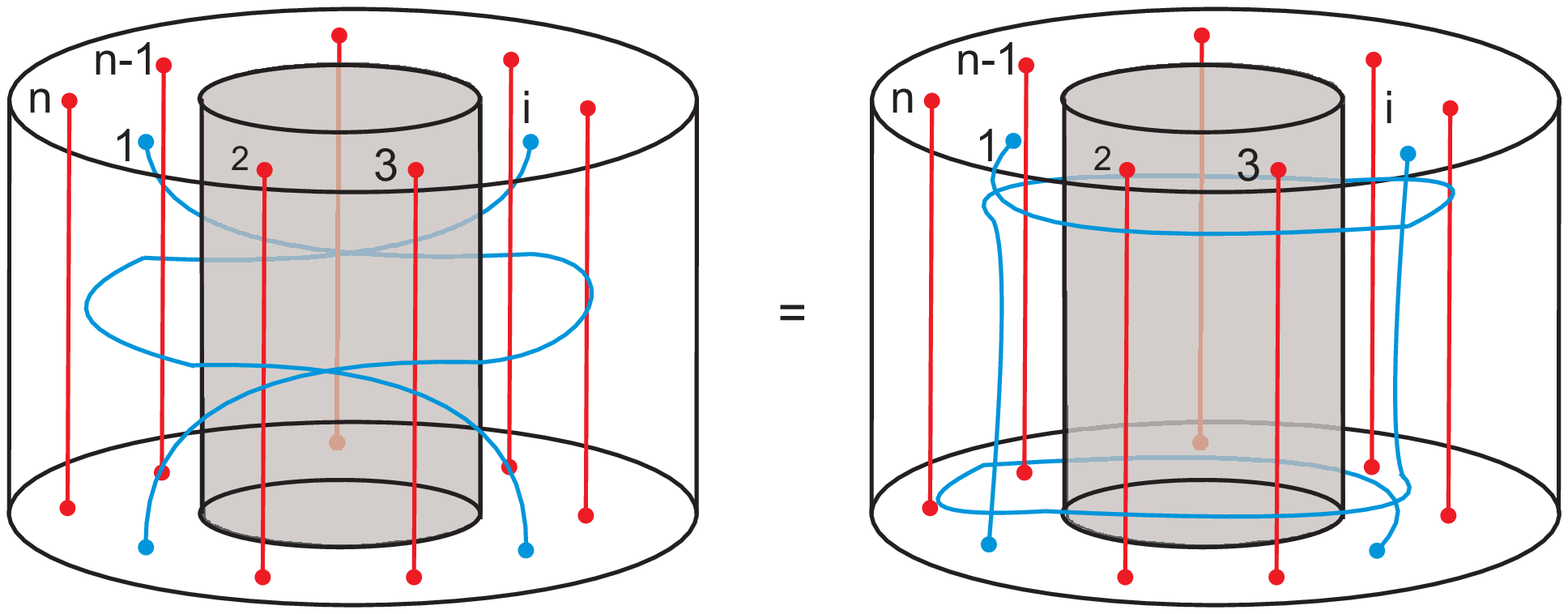}
\caption{The $i$-th generator of the group $F_n^0/R$ represented by
a balanced annular braid ($i=2,\dots ,n$).}\label{Fig_braid}
\end{figure}

\newpage

\begin{thm}\label{Thm_main}
\begin{itemize}
    \item The  group $F_n^0/R$   acts on the moduli spaces $ M_{2,3}( L).$

 \item The  group $F_n^0/R$ is isomorphic to the described below
     factor group of the
    \textit{
group of pure balanced annular braids}.
 Thus the elements of the  group $F_n^0/R$ can be represented by
 balanced annular braids. For instance, Fig. \ref{Fig_braid} depicts the generators
 of the group.
\item The  group $F_n^0/R$ is
isomorphic to $\mathbb{Z}^{n-1}$,
    and  is therefore commutative.

\end{itemize}
\end{thm}

Proof.

The first statement follows from the above discussion. To prove the
third statement, we construct an explicit homomorphism
$$\phi \colon F_n^0/R \to \Z^{n-1} \cong \{(w_1,w_2,\dots,w_n) \in \Z^n : \sum_{i=1}^n w_i=0.\}$$
 We start with
 the  \textit{balanced}
\textit{annular braid group} which is defined as follows:

$$B_n = \langle \Sigma_1, \Sigma_2, \dots \Sigma_n \mid \Sigma_i\Sigma_j = \Sigma_j\Sigma_i,$$$$
 \Sigma_i\Sigma_{i'}\Sigma_i = \Sigma_{i'}\Sigma_i\Sigma_{i'}
 \
\hbox{whereas}  \ i-j\neq \pm 1,\, i-i'= \pm 1 \rangle.$$ Next, we
take   the group $B_n^0$ of \textit{pure braids}, that is, the
kernel of the natural map $B_n \to S_n \hbox{  which maps } \Sigma_i
\hbox{ to } \sigma_i.$

As usual, we visualize a  braid as $n$ non-intersecting strands
living
 in a "thick" cylinder and going from the top to the bottom, see
 Fig. \ref{Fig_braid}.

 Finally, we introduce the group  $\overline{B_n^0}$, that is, the group $B_n^0$ factorized
 by all relations of type $(\Sigma_i)^2=1$.
The factorization means that the strands can pass freely through
each other, but not through the central part of the cylinder.

There is a natural isomorphism
$$ \psi: F_n^0/R \rightarrow \overline{B_n^0}$$
which maps $ s_i$ to $  \Sigma_i$.

 Besides, there is a  homomorphism
$$w\colon \overline{B_n^0} \to \Z^n,\, b \mapsto w(b) = (w_1(b), w_2(b), \dots, w_n(b))$$
where $w_i(b)$ is a winding number of the $i$-th strut of the braid
$b$ around the central part of the cylinder. It is easy to check
that for any pure braid $b$, we have
$$\sum_{i=1}^n w_i(b)=0.$$

Taken together, the two maps give the homomorphism
$$w \circ \psi :F_n^0/R \rightarrow \{(w_1,w_2,\dots,w_n) \in \Z^n : \sum_{i=1}^n w_i=0.\},$$
which is obviously bijective.

 Figure \ref{Fig_braid} depicts the preimage  of the vector
\newline $(1,0,0, \dots ,0,0,-1,0,0,\dots ,0)$ with just two non-zero
entries. The preimage of the vector in the group $F_n^0/R$ is
represented by $$s_{i+1}s_{i+2}\dots s_{i-1}s_{n-1}s_{n-2} \dots
s_2s_1.\qed$$

\comment{
\begin{cor}Denote by $SW_{n, \overline{i_0}}$ the subgroup of $SW_n$
generated by all $s_i$ except for  $s_1$. Then $SW_{n,
\overline{i_0}} \cong$
\end{cor}
}

\begin{prop}The critical points of the  function $A$ and $S$ (such that $S\neq 0$) are stable under the action of $F_n^0$.
\end{prop}

Proof.

 (2D). Critical points of the function $A$ are known to be cyclic polygons
 (see Theorem \ref{Thm_crirical_are_cyclic}). A cyclic
 polygon
 $P$ is completely determined by
 $r(P)$, $L$  and $E(P)$. The action of $F^0_n$ preserves them all.

   (3D). Assume that $P$ is a critical point such that $S(P)\neq 0$. Fix a polygon $P$ and a plane $\overrightarrow{S}^{\perp}$.
     First observe that a critical point is uniquely
    determined by radius $r(pr P)$ of the circumscribing circle,  the edge orientations $E(pr P)$,
    and the heights $h_i=dist(p_i, \overrightarrow{S}^{\perp}), \
    i=1,\dots , n$.

     Let $g$ be an element of $F_n^0/R$. Theorem \ref{Thm_crirical_3D} implies that
     the
    swap $s_i$ permutes the height differences $h_{i+1}-h_i$ and
    $h_i-h_{i-1}$.Therefore, $g$ maintains the height differences
    $h_{i+1}-h_i$. Besides, $g$ maintains both $E(pr P)$ and $r(pr
    P)$. By the above observation, $g$ maps $P$ to itself.\qed

\bigskip

Computer experiments show the following:

\def\scale{1}

\newcounter{no2}
\setcounter{no2}{0}

\def\Pic
{
\begin{minipage}[b]{0.12\linewidth}
\centering
\input p2\arabic{no2}.tex
\addtocounter{no2}{1}%
\end{minipage}
}

\centering
{
\begin{figure}
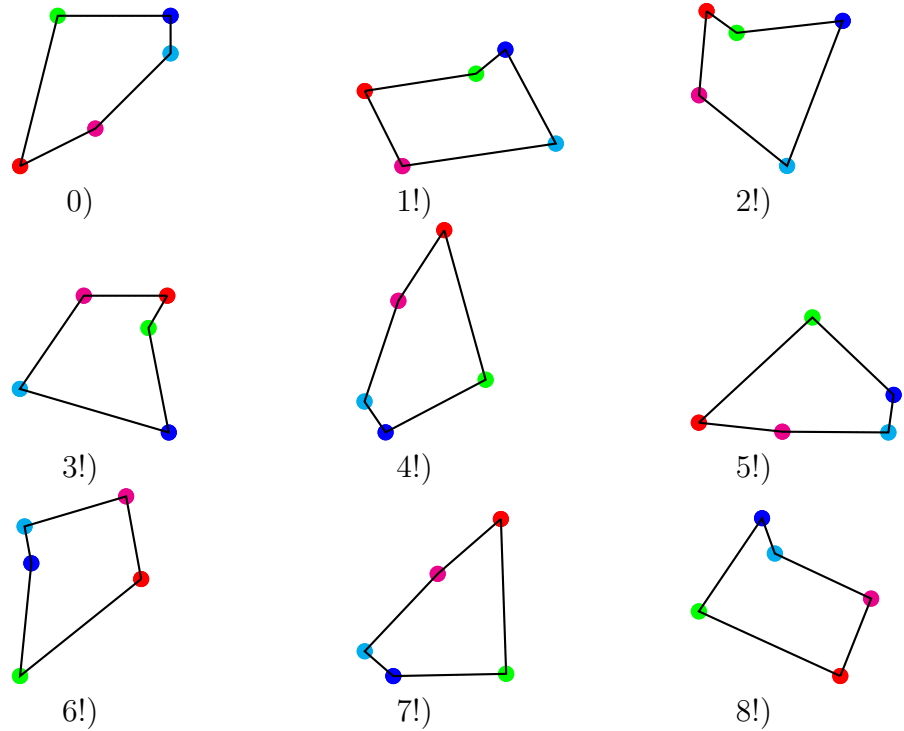

\begin{tabular}{ccc}
\hskip1cm\,\Pic\hskip1cm\,  &\hskip1cm\, \Pic\hskip1cm\, & \hskip1cm\,\Pic\hskip1cm\,\\\
0)& 1!)& 2!)\\
\hskip1cm\,\Pic\hskip1cm\,  &\hskip1cm\, \Pic\hskip1cm\, & \hskip1cm\,\Pic\hskip1cm\,\\\
3!)& 4!)& 5!)\\
\hskip1cm\,\Pic\hskip1cm\,  &\hskip1cm\, \Pic\hskip1cm\, & \hskip1cm\,\Pic\hskip1cm\,\\\
6!)& 7!)& 8!)\\
\end{tabular}
\caption{ We depict here the polygons (0) $P$ and the iterated
actions $g^{k!}(P)$} \label{cbgr}
\end{figure}}

\bigskip

\begin{Ex}
Let $g=s_4s_3s_2s_1$. There exists a pentagon $P$  such that $g^kP$
are all different for $k=1,2,\cdots, 8!$ (see Fig. 5).
\end{Ex}

\bigskip

{
\def\sr{{\color{red} s_1}}
\def\sg{{\color{green} s_2}}
\def\sb{{\color{blue} s_3}}
\def\sc{{\color{cyan} s_4}}
\def\sm{{\color{magenta} s_5}}

\def\scale{1}

\newcounter{no}
\setcounter{no}{9}

\def\Pic
{
\begin{minipage}[b]{0.12\linewidth}
\centering
\input \arabic{no}.tex
\addtocounter{no}{1}%
\end{minipage}
}

\centering
{
\begin{figure}
\begin{tabular}{ccc}
\hskip1cm\,\Pic\hskip1cm\,  &\hskip1cm\, \Pic\hskip1cm\, & \hskip1cm\,\Pic\hskip1cm\,\\\
a)& b)& c)\\
\hskip1cm\,\Pic\hskip1cm\,  &\hskip1cm\, \Pic\hskip1cm\, & \hskip1cm\,\Pic\hskip1cm\,\\\
d)& e)& f)\\
\end{tabular}
\caption{ The action of the first generator. We depict here a)
$P$,\quad b) $P \to \sr(P)$,\quad c) $\sr(P) \to \sg\sr(P)$,\quad
\dots, \quad e) $\sb\sg\sr(P) \to \sc\sb\sg\sr(P)$, \quad f)
$\sc\sb\sg\sr(P)$.} \label{cbgr}
\end{figure}}
\setcounter{no}{499}
\centering
{
\begin{figure}
\begin{tabular}{ccc}
\hskip1cm\,\Pic\hskip1cm\,  &\hskip1cm\, \Pic\hskip1cm\, & \hskip1cm\,\Pic\hskip1cm\,\\\
a)& b')& c')\\
\hskip1cm\,\Pic\hskip1cm\,  &\hskip1cm\, \Pic\hskip1cm\, & \hskip1cm\,\Pic\hskip1cm\,\\\
d')& e')& f')\\
\end{tabular}
\caption{The action of the second generator.  We depict a) $P$,\quad
b') $P \to \sg(P)$,\quad c') $\sr(P) \to \sb\sg(P)$,\quad \dots,
\quad e') $\sc\sb\sg(P) \to \sm\sc\sb\sg(P)$, \quad f')
$\sm\sc\sb\sg(P)$.} \label{mcbg}
\end{figure}}

\begin{Ex}

For the pentagon $P$ depicted in Fig.  6, the pentagons \
$\sc\sb\sg\sr(P)$ and $\sm\sc\sb\sg(P)$ are different. This means
that they are different for a
 generic pentagon.
\end{Ex}
}
 These two examples motivate the following conjecture:

\begin{conj} For a generic polygonal linkage, the groups $Stab $ and $R $ coincide,  i.e., $SW_n = F^0_n / R= \mathbb{Z}^{n-1}$.
\end{conj}


\end{document}